\magnification=1200
\baselineskip=16pt plus 2pt minus 1pt

\font\magbf=cmbx10 scaled\magstep2


\def\mbox#1{\leavevmode\hbox{#1}}
\def\to{\mbox{$\longrightarrow$}}	
\def\frac#1#2{{#1\over #2}} 
 

\def\zz{\mbox{${\bf Z}$}}	
\def\qq{\mbox{${\bf Q}$}}	
\def\calc{\mbox{${\cal C}$}}
\def\calg{\mbox{${\cal G}$}} 

\count35=0
\count36=0
\count38=0


\def\section#1{\advance\count35 by 1
\vskip4ex\noindent{\magbf\number\count35 \ \ \ #1}
\count36=0\count38=0 \hfill\vskip2ex}

\def\subsection{\vskip1ex\advance\count36 by 1
{\noindent{\bf\number\count35.\number\count36}\ }} 

\def\today{\ifcase\month\or January\or February\or March\or  April\or
May\or
June\or July\or August\or September\or October\or November\or December\fi
\space\number\day, \number\year}

\def\beginth#1{\vskip3ex\advance\count36 by 1 {\noindent{\bf
\number\count35 .\number\count36} \ \ {\bf #1.} }\begingroup\sl} 

\def\endth{\endgroup\hfill\vskip2ex}

\def\defn{\vskip3ex\advance\count36 by 1 {\noindent{\bf \number\count35
.\number\count36} \ \ {\bf Definition.} }}

\centerline{\magbf   Knot Concordance and Torsion}

\vskip3ex
\centerline{Charles Livingston and Swatee Naik}

\centerline{November 30, 1999}

\vskip7ex
\section{Introduction}

The classical knot concordance group, $\calc_1$, was defined in 1961 by
Fox [F].  He  proved that it is nontrivial by
finding elements of order two;  details were presented in [FM].  Since
then one of the most vexing questions concerning the concordance group
has been whether it contains elements of finite order other than
2--torsion.  Interest in this question was heightened by Levine's proof
[L1, L2] that in all higher odd dimensions the knot concordance group
contains an infinite summand generated by elements of order 4.  In our
earlier work studying this problem we proved the following [LN]:

\beginth{Theorem}   Let $K$ be a knot in $S^3$ with 2--fold branched cover
$M_K$.  If  the order of the first homology with integer coefficients
satisfies $|H_1(M_K)| = pm$ with
$p$ a prime congruent to 3  mod 4 and gcd(p, m) = 1, then $K$ is of
infinite order in the classical knot concordance group, 
$\calc_1$.\endth 

\noindent
An immediate corollary was that all of the prime knots of
less than 11 crossings that are of order 4 in the algebraic concordance
group are of infinite order in the concordance group. There are 11 such
knots [M]. One simple case of a much deeper corollary states that if the
Alexander polynomial of a knot satisfies 
$\Delta_K(t) = 5t^2 - 11t + 5$  then $K$ is of infinite order in
$\calc_1$.  According to Levine [L2], any higher dimensional knot with this
polynomial is of order 4 in concordance. 

\vskip3ex
Here our goal is to prove the following enhancement of the theorem
stated above:

\beginth{Theorem}   Let $K$ be a knot in $S^3$ with 2--fold branched cover
$M_K$.  If $H_1(M_K) = \zz_{p^n} \oplus G$ with $p$ a prime congruent to
3  mod 4, $n$ odd, and $p$ not dividing the order of $G$, then $K$ is of
infinite order in  $\calc_1$.\endth

\noindent
As we will describe below, the significance of this result goes beyond
its apparent technical merit; however, even in terms of computations it
is an important improvement. 
Let $H_p$ denote the $p$--primary summand of $H_1(M_K)$.  

\beginth{Corollary}Let $n$ be a positive integer such that
some prime $p$ congruent to 3 mod 4 has odd
exponent in the prime power factorization of
$4n+1$. Then a
knot $K$ with Alexander polynomial $nt^2 -(2n+1)t +n$ and $H_p$
cyclic  is of infinite order in the 
concordance group.\endth

\noindent
Note that according to Levine [L2], any such knot represents an order 4
class in the algebraic concordance group.
The $n$--twisted doubles of 
knots provide infinitely many examples of knots
with Alexander polynomial $nt^2 -(2n+1)t +n$ and $H_1(M_K)$ cyclic. 
Further details and examples will be provided in the last section.

Casson and Gordon's first examples of algebraically slice knots that
are not slice [CG1, CG2] were taken from the set of
twisted doubles of the unknot.  
Our analysis extends theirs to a
much larger class of knots and is not restricted to doubles of the
unknot.  Notice here the rather remarkable fact that an abelian
invariant of a knot is being used to obstruct an algebraically slice
knot from being slice!

Theorem 1.2 is  relevant to deeper
questions concerning the concordance group.  The underlying conjecture
is that the only torsion in the knot concordance group is 2--torsion,
arising from amphicheiral knots (see related questions in [G, K1, K2]); a
positive solution to this conjecture seems far beyond the tools now
available to study concordance, in either the smooth or topological
locally flat category.  However, two weaker conjectures are possible.

\beginth{Conjecture} If a knot $K$  represents 4--torsion in the
algebraic concordance group, then it is of infinite order in
concordance.  \endth

\noindent
A simpler conjecture is:

\noindent
\beginth{Conjecture}  There exists a class of order 4 in the
algebraic concordance group that cannot be represented by a knot of order
4.  In particular, Levine's homomorphism does not split. \endth 

Theorem 1.1 provided candidates for verifying Conjecture 1.5 but there
are two difficult steps to extending that result from a representative of
an algebraic concordance class to the entire class.  
It is a consequence of Witt theory (that we won't
be using elsewhere in this paper) that such an extension will have two
parts: one must be able to handle the case where $H_p = \zz
_{p^n}$, with $n > 1$,  and also the case where $H_p$ is a direct
sum of such factors.

The results of this
paper deal with 
the first part of the extension problem.  A number of special
cases of direct sums have been successfully addressed
by us in unpublished work, but
the necessary general result for sums has not yet been achieved.  When
it is, that result along with Theorem 1.2 should provide a proof of
Conjecture 1.5 and perhaps 1.4 .

The work of this paper is largely algebraic.  In the next section we
will summarize the topological results that we will be using.  All the
work that appears here applies in both the topological locally flat and
the smooth category.  
In Section 3 we give a proof of Theorem 1.2. The proof is fairly
technical  and extends the techniques used in proving Theorem 1.1 in
[LN].  Section 4 discusses examples.

\section{Background and notation}
\subsection{\bf Knots and the concordance group}

We will work in the smooth category, but as just mentioned, all results
carry over to the topological locally flat setting.  Homology groups
will always be with $\zz$ coefficients unless otherwise mentioned.

A knot is formally defined to be a smooth oriented pair,
$(S^3,K)$, with $K$ diffeomorphic to $S^1$.  We will denote such a pair
simply by $K$. 
A knot $K$ is called {\sl slice} if $(S^3, K) = \partial (B^4, D)$, where
$D$ is a smooth 2--disk properly embedded in the 4--ball $B^4$.
Knots $K_1$ and $K_2$ are called {\sl
concordant} if $K_1\ \#\ \mbox{--}K_2$ is slice, where $\mbox{--}K$
represents the mirror image of
$K$, formally $(\mbox{--}S^3,\mbox{--}K)$. The set of concordance classes
of knots forms an abelian group under connected sum, denoted $\calc_1$. 
The order of $K$ in the knot concordance group is hence the least
positive
$n$ for which the connected sum of $n$ copies of $K$, $nK$, is slice.

Levine defined a homomorphism of $\calc_1$ onto a group, $\calg$, that
is isomorphic to the infinite direct sum, $\calg \cong \zz^\infty
\oplus\zz_2^\infty
\oplus\zz_4^\infty$.  For higher dimensions the corresponding
homomorphism is an isomorphism, but in dimension 3 there is a nontrivial
kernel, as first proved by Casson and Gordon [CG1, CG2].  For details
concerning $\calg$, the {\sl algebraic concordance group}, see [L1, L2].

\subsection{\bf Casson-Gordon invariants and linking forms}
	
Let $M_K$ denote the 2--fold branched cover of $S^3$ branched over $K$,
and let $\chi$ denote a homomorphism from 
$H_1(M_K)$ to 
$\zz_{p^k}$ for some prime $p$.  The Casson-Gordon invariant,
$\sigma(K,\chi)$ is then a well defined rational invariant of the pair
$(K,\chi)$.  (In Casson and Gordon's original paper, [CG1], this
invariant is denoted $\sigma_1 \tau(K,\chi)$, and $\sigma$ is used for
a closely related invariant.)

On any rational homology sphere, such as $M_K$, there is a nonsingular
symmetric linking form, $\beta :H_1(M_K) \rightarrow \qq / \zz$. As
before, let
$H_p$ be the $p$--primary summand of $H_1(M_K)$.  
The main result
in [CG1]
concerning Casson-Gordon invariants and slice knots that we will be
using is the following:

\beginth{Theorem} If $K$ is slice there is a subgroup (or {\it
metabolizer}) 
$L_p \subset H_p$ with $|L_p|^2 = |H_p|$, \ $\beta(L_p,L_p) = 0$, and
$\sigma(K,\chi) = 0$ for all
$\chi$ vanishing on $L_p$.\endth

\noindent
We will also need the additivity theorem proved by Gilmer [Gi].

\beginth{Theorem}  If $\chi_1$ and $\chi_2$ are defined on $M_{K_1}$
and   $M_{K_2}$, respectively, then we have
$\sigma( K_1\ \#\ K_2, \chi_1\ \oplus \
\chi_2) =
\sigma(K_1,\chi_1) + \sigma(K_2,\chi_2)$. \endth

\noindent
Any homomorphism $\chi$ from $H_p$ to $\zz_{p^r}$ is
given by linking with some element $x \in H_p$.  
In this situation we have the following (see Section 4 of [LN]).

\beginth{Theorem}  If $\chi \colon H_p \to \zz_{p^r}$
is a character obtained by linking with the element $x \in H_p$,
then $\sigma(K,\chi) \equiv \beta(x,x)\ \mbox{\rm modulo}\  \zz$.   \endth

\noindent
A simple corollary, using the nonsingularity of $\beta$ is:

\beginth{Corollary} If $H_p = \zz_{p^n}$ and  $\chi$ maps onto 
$\zz_{p^k}$ with
$k > n/2$ then $\sigma(K,\chi) \ne 0$. \endth

\noindent
Finally, we will use the result below which is a consequence of the
fact that the linking form $\beta$ 
gives a map from $H_p$
onto Hom$(L_p, {\bf Q}/{\bf Z}) \cong L_p$, with kernel equal to $L_p$.

\beginth{Theorem} 
With $H_p$ and $L_p$ as in Theorem 2.3, we have
$H_p / L_p \cong L_p .$ \endth

\section{Proof of Theorem 1.2}

Let $K$ be a knot in $S^3$ with the 2--fold branched cover $M_K$.
Suppose that $H_1(M_K) = \zz_{p^n} \oplus G$ with $p$ a prime congruent to
3  mod 4, $n$ odd, and $p$ not dividing the order of $G$.
We want to show that $K$ is of
infinite order in  $\calc_1$.
The linking form of
$H_1(M_K)$ represents an element of order 4 in the Witt group of $\zz _p$
linking forms. (See Corollary 23 (c) in [L2].)
If $K$ is of concordance order $d$, since Levine's homomorphism maps the
concordance class of $K$ to an order 4 class, we have $d = 4k$, for some
positive integer $k$. We must analyze the possible metabolizers
$L_p$ for
$(\zz_{p^n})^{4k}$.

A vector in $L_p$ can be written as
$x = (x_i)_{i=1\ldots d} \in (\zz_{p^n})^d$. 
Applying the Gauss-Jordan algorithm to a basis for $L_p$,
and perhaps reordering, we can find a generating set of a particularly
simple form.  The next example illustrates a possible form for one such
generating set, where the generators appear as the rows of the matrix.

\beginth{Example}
Let $H_p = \left( \zz _{p^3}\right) ^8$. A generating
set for some metabolizer $L_p$ of the standard nonsingular
$\qq / \zz$ linking form can be written as follows:\endth
\vskip-.1in
$$ \left( \matrix{
1 & * & * & * & * & * & * & * \cr
0 & p & 0 & 0 & * & * & * & * \cr
0 & 0 & p & 0 & * & * & * & * \cr
0 & 0 & 0 & p & * & * & * & * \cr
0 & 0 & 0 & 0 & p^2 & 0 & 0 & * \cr
0 & 0 & 0 & 0 & 0 & p^2 & 0 & * \cr
0 & 0 & 0 & 0 & 0 & 0 & p^2 & * \cr
         } \right) $$
\vskip.1in

\noindent In the above matrix, there is 1 row corresponding to $p^0$, 
and 3 rows each for $p^1$ and $p^2$. 

We will denote the
number of rows corresponding to $p^i$ by $k_i$, the 
vectors in these $k_i$ rows by $v_{i,1}, \cdots , v_{i, k_i}$,
and $\sum_{j=0}^i k_j$ by $S_i$.  For notational purposes, let $S_{-1} = 0$. 
Then, in general, the generating set consists of 
$\{v_{i,j} \}_{i=0, \ldots ,n-1, j =1, \ldots ,k_i}$ where 
$0 \le k_i \le 2k$, 
such that the first $S_i$ 
entries of 
$v_{i,j}$ are 0, except for the 
$ S_{i-1}
+j$ entry which is $p^i$,
and each of the remaining entries is divisible by $p^i$.
From 2.7 it follows that $k_i = k_{n-i},$ for $i > 0,$ and
$S_{(n-1)/2}
= 2k$.

\beginth{Definition} If $a \in H_p$, let $\chi_a \!: H_1(M_K)
\rightarrow  \qq / \zz$ be the character given by
linking with $a$.  In the case that $H_p $ is cyclic, isomorphic
to $\zz_{p^n}$, we can fix a generator of $H_p$ and write  $\chi_a$
where $a$ is an integer representing an element in $\zz_{p^n}$.  
\endth

\noindent
With this
notation, we now see that our goal is to show that
$\sigma(K,\chi_{p^{(n-1)/2}})= 0$. Since $\chi_{p^{(n-1)/2}}$
maps onto $\zz_{p^{(n+1)/2}}$ this will contradict 2.6
and it will follow that $K$ cannot be of finite order.

As in Example 3.1, arrange the $\{v_{i,j} \}$ as rows of a 
$(4k-k_0) \times 4k $ matrix following the dictionary 
order on $(i,j)$.
We multiply the first $k_0$ vectors by $p^{n-1}$,
 the next $k_1$  vectors by $p^{n-2}$,
and, so on, to obtain 
$p^{n-1}$ on the diagonal. Clear the off-diagonal entries in the left
$(4k-k_0) \times (4k-k_0)$ block. Now, adding
all the rows gives us a vector in $ L_p$ with the first
$4k-k_0$ entries equal to $p^{n-1}$. This vector corresponds 
to a character $\chi ,$ 
given by linking an element with it, to $\zz _p$
on which the Casson-Gordon
invariants should vanish. That is, 
$\sigma ((4k)K,\chi) = 0$. By 2.4
this leads to a relation of the form
$(4k-k_0) \sigma (K,\chi_{p^{n-1}} )+ 
\sum_{x_i \ne 0} \sigma (K,\chi_{x_i}) = 0$, where $x_i$ are 
the remaining $k_0$ entries, each of which is divisible
by $p^{n-1}$.

The set of nonzero characters from $\zz _{p^n}$ to $\zz _p$ is isomorphic
to the multiplicative group of units in $\zz _p$, which is a cyclic
group of order $p-1$. 
Denoting a generator for this group by $g$,
each nonzero $\chi _{x_i}$ corresponds to
$g^{\alpha_i}$ for some $\alpha_i$. 
The correspondence can be given by $\chi_{x_i} \to g^{x_i / p^{n-1}}$.
As in [LN] we use further shorthand,
setting $t^{\alpha_i} = \sigma(K,\chi_{x_i})$.
Each metabolizing vector leads to a relation $\sum_{x_i \ne 0} t^{\alpha_i} = 
0$.  
Note that at this point the symbol $t^{\alpha}$ does not represent a
power of any element ``$t$'', it is purely symbolic.  However it does
permit us to view
 the relations as being elements in the ring
$\zz[\zz_{p-1}]$.  Furthermore, since $\sigma(K,\chi_{x_i})
 = \sigma(K,\chi_{p^n -x_i})$, we have
that
$t^j = t^{j + (p-1) /2}$.  (Recall that $g^{(p-1)/{2}} = -1$.)  Hence, we
can view the relations as sitting in $\zz[\zz_{q}]$, where $q = (p-1)/2$.

If a metabolizing vector $x$ corresponds to the relation
$f = 0$, where $f$ is represented by an element in $\zz[\zz_q]$,  then 
$ax$ corresponds to the relation $t^\alpha f$ where $g^\alpha = a$. 
It follows that the relations between Casson-Gordon invariants 
generated by  the 
element $x \in L_p$ together with 
its multiples form an ideal in $\zz[\zz_q]$ generated by the
polynomial $f$.

With this in mind our relation can be written as
$f = (4k - k_0) + \sum_{i=1}^{k'}t^{\alpha_i} = 0$, where
$k' \le k_0 $. (Note that $4k - k_0 = S_{n-1}.$) 
We 
show that the ideal generated by
$f$ in
$\zz[\zz_q]$ contains a nonzero integer.  This will follow from the fact  
that 
$f$ and $t^q -1$ are relatively prime, which will be the case unless $f$ 
vanishes at some 
$q$--th root of unity, say 
$\omega$; however, by considering norms and the triangle inequality we see
that this will be the case only if $k' = 2k$ and
$\omega^{a_i} = $ --1 for all $i$.  But since $q$ is odd, no power of
$\omega$ can equal --1.

It follows that $n \sigma(K,\chi_{p^{n-1}})= 0$, 
for some $n \in \zz$.
This implies that
$\sigma(K,\chi_{p^{n-1}})= 0$. Similarly we can show that
$\sigma(K,\chi_{ap^{n-1}})= 0$, for $0 < a < p$.

Next, let $l$ be a nonnegetive integer, and
assume that $\sigma(K,\chi_{ap^s})= 0$, for all $a \in \zz$, and
all $s$ such that $l < s \le n-1$.
We will show that $\sigma(K,\chi_{p^l})= 0$.

For $ 0 \le i \le S_l$,
we multiply the vectors
from the $(S_{i-1}+1)$st to the $S_i$th vector by $p^{l 
-i}$, clear 
off-diagonal entries in the upper left 
$S_l \times S_l$ square block, and add the first $S_l$
rows to get
a vector in $L_p$  with first
$S_l$ entries equal to $p^l 
$, and the remaining entries divisible by $p^l 
$. Since we have assumed that $\sigma(K,\chi_{ap^s
})= 0$, for $l < s \le n-1$,
we can ignore the entries which are of the form $ap^s 
$, with $s > l$.
Then we have a character to the multiplicative group
of units in $\zz _{p^{n-l}
}$. Since $p$ is odd, this is a cyclic group
of order $p^{n-l-1} (p-1)$ (see [D]). Again, since
$\sigma(K,\chi_{x_i}) 
 = \sigma(K,\chi_{p^n -x_i})$, we 
can view the relations as sitting in $\zz[\zz_{q}]$, where $q = p^{n-l-1}
(p-1)/2$.
As $p^{n-l-1}
(p-1)/2$ is odd, as above, it follows that the relation 
$f = S_l
 + \sum_{i=1}^{k'}t^{\alpha_i} = 0$, where $
0 \le k' \le 4k - S_l 
$, is relatively prime to $t^q -1$. It follows that 
$\sigma(K,\chi_{p^l})= 0$.  Similarly, $\sigma(K,\chi_{a p^l})= 0$ for $0 <
a < p$.

Thus, we have 
$\sigma(K,\chi_{p^{(n-1)/2}})= 0$, which contradicts 2.6,
and proves that $K$ cannot be of finite order in the concordance
group.

\section{Examples}

Basic examples illustrating the applicability of Theorem 1.2 are easily
constructed.  For instance, since the 2--fold branched cover of $S^3$
over an unknotting number one knot has cyclic homology, to apply
Theorem 1.2 we only need to check the order of $H_1(M_K)$ which 
equals the Alexander polynomial evaluated at $-1$. We have the 
following. 

\beginth{Corollary} Let $K$ be an unknotting number one knot with 
Alexander polynomial $\Delta$. If a prime $p$ which is congruent
to 3 mod 4 appears in the prime power factorization of
$\Delta(-1)$ with an odd exponent, then $K$ is of 
infinite order in the concordance group. \endth 

\noindent 
More generally, suppose that there is a
3--ball $B \subset S^3$ intersecting the knot $K$ in  two arcs so that the
 pair $(B, B \cap K)$ is trivial and so that removing  $(B, B \cap K)$
from $S^3$ and gluing it back in via a homeomorphism of the boundary yields
the unknot.  Since the 2--fold branched cover of the ball over two trivial
arcs is a solid torus,  the 2--fold branched cover of
$S^3$ over $K$ is formed from $S^3$  (the 2--fold branched cover of $S^3$
over the unknot) by removing a solid torus and sewing it back in via some
homeomorphism. In particular, the 2--fold branched cover has cyclic
homology.  Such knots include all unknotting number one knots
and all 2--bridge knots. In the case of a 2--bridge knot $K(p,q)$, we have
$H_1(M_K) = \zz _p$. 

\beginth{Corollary} The 2--bridge knot $K(p,q)$ has infinite order in the
knot concordance group if some prime congruent to 3 mod 4 has odd
exponent in $p$.\endth

\noindent
The following theorem, an
immediate consequence of a result of Levine (Corollary 23 in [L2]), provides
us with more examples of knots which represent torsion in the
algebraic concordance group.
\vfill\eject
 
\beginth{Theorem}  If a knot $K$ has quadratic Alexander
polynomial
$\Delta(t)$ then:

\item{(a)} $K$ is of finite order in the algebraic concordance group if
and only if $\Delta(1)\Delta(-1) <0$, in which case $K$ is of order 1, 2
or 4.

\item{(b)} $K$ is of order 1 if and only if $\Delta(t)$ is reducible.

\item{(c)} if $K$ is finite order, and $\Delta(t)$ is irreducible, then
$K$ is of order 4 in the algebraic concordance group if and only if
for some prime $p >0$ with $p \equiv 3\ \mbox{\rm mod}\ 4$, 
$\Delta(1)\Delta(-1) = - p^a q$ where $a$ is odd and $q >0$ is
relatively prime to $p$. \endth   

\noindent
Consider the  knot $K_a$, the $a$--twisted double of
some knot $K$.
The Seifert form for this knot is 
$$ V =\left( \matrix{
a & 1\cr
0& -1\cr
} \right),$$
it has Alexander polynomial $\Delta(t) = at^2 - (2a+1)t +a$, and the
homology of the 2--fold cyclic branched cover is $\zz_{|4a+1|}$.  Levine's
result, Theorem 4.3, applies to determine the algebraic order of all of
these knots. 
(In the case that $K$ is unknotted, $K_a$ can be described as the
2--bridge knot $K(4a+1, 2a)$.) 

\beginth{Corollary} The $a$--twisted double of a knot $K$:

\item{(a)} is of infinite order in the algebraic concordance group, $\calg$,
if
$a < 0$.

\item{(b)} is algebraically slice if $a > 0$ and $4a+1$ is a perfect square.

\item{(c)} is of order 2 in $\calg$ if $a > 0$, $4a+1$ is not a perfect
square, and every prime congruent to 3 mod 4 has even exponent in the
prime power factorization of $4a+1$.

\item{(d)} is of order 4 if  $a >0$ and some prime congruent to 3 mod 4 has 
odd exponent in $4a+1$.

 \endth

\noindent
Casson and Gordon [CG1, CG2] proved that if $K$ is unknotted, then all knots
covered by case (b) above are actually of infinite order in concordance,
except if
$a = 2$ in which case $K_2$ is slice.  An immediate consequence of Theorem
1.2 is:

\beginth{Corollary} If $K_a$ is of order 4 in $\calg$ then it is of
infinite order in the knot concordance group. \endth

\noindent
As in Corollary 9.5 of [LN] a simple argument using Corollary 4.5
gives an infinitely generated free subgroup of ${\cal C} _1$
which consists of of algebraic slice knots.
(It was first shown by Jiang in [J]
that the kernel of Levine's homomorphism is infinitely generated.) 
The extensive calculations of [CG1, CG2] are here replaced with a
trivial homology calculation. 
Moreover, the results
of [CG1, CG2] apply only in the case that $K$ is unknotted, a restriction
that does not appear in Corollary 4.5. 

Recently, Tamulis [T] has proved that
in the case that $K$ is unknotted, if $K_a$ is of order 2 in $\calg$ and 
$4a+1$ is prime, then
$K_a$ is of infinite order in concordance. 
\vskip1ex
\noindent{\bf Counterexamples.\ }  Given these previous examples, it is a bit
unexpected that Theorem 1.2 does not apply in all cases of order 4 knots
with quadratic Alexander polynomial.  The difficulty is that the conditions
of Theorem 4.3 do not assure that the homology of the 2--fold cover is
cyclic. The next example demonstrates this.  It is the simplest possible
example in terms of the coefficients of the Alexander
polynomial; its complexity illustrates the strength of
Theorem 1.2.  The example is obtained by letting
$K$ be a knot  with Seifert form:

$$ V =\left( \matrix{
21 & 53\cr
52 & 21 \cr
} \right). $$

The Alexander polynomial for $K$ is $\Delta(t) = 2315 -4631t +2315t^2$.
We have that $\Delta(1) = -1$, $\Delta(-1) = 9261 = 3^3 7^3$, and hence
by Theorem 4.3, $K$ is of order 4 in the algebraic concordance group.

The homology of $M_K$ is presented by $V+V^t$:
$$ V =\left( \matrix{
42 & 105\cr
105 & 42 \cr
} \right). $$
A simple manipulation shows that this presents the group $\zz_3 \oplus
\zz_9 \oplus \zz_7 \oplus \zz_{49}$.  Because this is not cyclic, Theorem
1.2 does not apply. As mentioned in the introduction,
we have been able to extend our results to special cases of direct sums of
cyclic groups, and one of those extensions applies to the group 
$\zz_3 \oplus \zz_9$. Hence it can actually be shown that any knot with
this Seifert form is not of order 4 in concordance.

\vfill
\eject

\vskip8ex
\centerline{\magbf References}
\vskip4ex

\item{[CG1]} A. Casson and C. Gordon, {\sl Cobordism of classical knots,}
Preprint, Orsay (1975). (Reprinted in ``A la recherche de la Topologie perdue'', ed. by 
Guillou and Marin,
Progress in Mathematics, Volume 62, Birkhauser, 1986.) 

\item{[CG2]} A. Casson and C. Gordon, {\sl On slice knots in dimension
three}, in Proc. Symp. Pure Math. 32 (1978), 39--54.

\item{[D]} D. Dummit and R. Foote, {\sl Abstract Algebra, 2nd ed.},
Prentice Hall, New Jersey, 1999.

\item{[F]} R. Fox, {\sl A quick trip through knot theory}, 1962 Topology
of 3--manifolds and related topics (Proc. The Univ. of Georgia Institute,
1961)  120--167 Prentice-Hall, Englewood Cliffs, N.J.

\item{[FM]} R. Fox and J. Milnor, {\sl Singularities of $2$--spheres in
$4$--space and cobordism of knots}, Osaka J. Math. 3 (1966) 257--267.

\item{[G]} C. Gordon, {\sl Problems}, in Knot Theory, ed. J.-C. Hausmann,
Springer Lecture Notes no. 685, 1977.

\item{[Gi]} P. Gilmer,  {\sl Slice knots in $S\sp{3}$},  Quart. J. Math.  
Oxford Ser. (2) 34 (1983), no. 135, 305--322.

\item{[J]} B. Jiang, {\sl A simple proof that the concordance group of
algebraically slice knots is infinitely generated}, Proc. Amer. Math. Soc.
83 (1981), 189--192.

\item{[K1]} R. Kirby, {\sl Problems in low dimensional manifold theory},
in Algebraic and Geometric Topology (Stanford, 1976), vol 32, part II of
Proc. Sympos. Pure Math., 273--312. 

\item{[K2]} R. Kirby, {\sl Problems in low dimensional manifold theory}, 
Geometric Topology, AMS/IP Studies in Advanced Mathematics, ed. W. Kazez,
1997.

\item {[L1]} J. Levine, {\sl Knot cobordism groups in codimension two},
Comment. Math. Helv. 44 (1969), 229--244.

 \item {[L2]} J. Levine, {\sl Invariants of knot cobordism}, Invent.
Math. 8 (1969), 98--110.

\item  {[LN]} C. Livingston and S. Naik, {\sl Obstructions to 
4--Torsion in the Classical Knot Concordance Group},  J. Diff. Geom.

\item{[M]} T. Morita, {\sl Orders of knots in the algebraic knot cobordism
group}, Osaka J. Math. 25 (1988), no. 4, 859--864.

\item{[R]} D. Rolfsen, {\sl Knots and Links}, Publish or Perish 7,
Berkeley CA (1976).

\item{[T]} A. Tamulis, {\sl Concordance of classical knots}, thesis, Indiana
University Department of Mathematics, 1999.

\end